\documentclass[11pt,letterpaper]{article}
\usepackage[letterpaper, margin=0.75in]{geometry}

\usepackage{graphicx,color}
\usepackage[cp850]{inputenc}
\usepackage[T1]{fontenc}
\usepackage{rotating}
\usepackage[f]{esvect}
\usepackage{amsmath, amsthm, amssymb}
\usepackage{rotating}
\usepackage{lscape}
\usepackage{mathrsfs}
\usepackage{subcaption}
\usepackage{mathtools}
\usepackage[toc,page]{appendix}
\usepackage{hyphenat}
\usepackage{mathptmx}
\usepackage[scaled=.65]{helvet}
\usepackage{courier}
\usepackage[right, mathlines]{lineno}

\usepackage[normalem]{ulem}

\DeclarePairedDelimiter{\ceil}{\lceil}{\rceil}

\def\tablenotes{\bgroup\parfillskip=0pt plus 1fil
\leftskip=0pt\relax \rightskip=0pt
\vskip2pt\footnotesize}
\def\endtablenotes{\vskip1pt\egroup}

\newtheorem{theorem}{Theorem}[section]

\renewcommand{\epsilon}{\varepsilon}
\renewcommand{\leq}{\leqslant}
\renewcommand{\geq}{\geqslant}
\renewcommand{\d}{\mathrm{d}}

\renewcommand{\epsilon}{\varepsilon}

\newcommand{\TimeDeriv}{\frac{\textrm{d}}{\textrm{dt}}}

\allowdisplaybreaks

\usepackage[square,numbers,sort]{natbib}

\usepackage{authblk}
\usepackage{graphicx}
\usepackage{overpic}
\usepackage{tikz}
\usepackage{pgfplots}
\usepgfplotslibrary{fillbetween}
\pgfplotsset{compat=1.10}
\usepackage{setspace}
 \usepackage[right,mathlines]{lineno}
 
\let\oldequation\equation
\let\oldendequation\endequation

\renewenvironment{equation}
  {\linenomathNonumbers\oldequation}
  {\oldendequation\endlinenomath}
  
\let\oldalign\align
\let\oldendalign\endalign

\renewenvironment{align}
  {\linenomathNonumbers\oldalign}
  {\oldendalign\endlinenomath}

\renewenvironment{align*}
  {\linenomathNonumbers\oldalign\notag}
  {\notag \oldendalign \endlinenomath}

 \graphicspath{ {./Figures/} }

\begin{document}

\title{Genuine and spurious bistability in a simple epidemic model with waning immunity}
\author[1,$\sharp$,*]{Francesca Scarabel}  \author[1]{Harry Coldwell} \author[1,$\sharp$,*]{Tyler Cassidy} 
\affil[1]{School of Mathematics, University of Leeds, Leeds, LS2 9JT, United Kingdom}
\affil[ $\sharp$ ]{Equal contribution}
\affil[*]{Correspondence to: f.scarabel@leeds.ac.uk, t.cassidy1@leeds.ac.uk}
\date{}
\maketitle

\begin{abstract}
We study an infection-age structured epidemic model in which both the infectivity and the rate of loss of immunity depend on the time-since-infection. The model can be equivalently viewed as a nonlinear renewal equation for the incidence of infection or as a partial differential equation for the density of infected individuals. We explicitly consider gamma, rather than Erlang, distributed durations of infection using a combination of ODE approximations and numerical bifurcation methods. We show that the shape of this distribution strongly influences stability of the endemic equilibrium, even when the basic reproduction number $R_0$ and the mean duration of infectiousness are fixed. Moreover, we establish the existence of regions of bistability, where a stable endemic equilibrium coexists with a stable periodic orbit. To our knowledge, this provides the first example of bistability in infection-age structured models with waning immunity alone. Finally, we show how common compartmental modelling approaches, which impose implicit assumptions on the distribution of the duration of infection, can lead to spurious dynamical outcomes. Taken together, our analysis underscores the crucial role of distributional structure in epidemic modelling and provides new insights into the rich dynamics of infection-age structured SIS/SIRS models.

\end{abstract}

\clearpage

\section{Introduction}

Mathematical models have been extensively employed to study the dynamics of infectious diseases at the population scale over the past century \citep{Kermack1927}. The most common approach involves compartmental models formulated as systems of ordinary differential equations (ODEs) that track the time evolution of susceptible, infected, and recovered individuals \citep{Champredon2018,keeling2008modeling,Wearing2005,Krylova2013}. These models are straightforward to analyze and simulate, and provide useful tools to assess interventions such as vaccination or isolation \citep{Magal2010,Magpantay2016,Magpantay2017,Le2021impact}. However, the classical ODE framework often assumes exponentially or Erlang distributed sojourn times and constant infectivity, which may not accurately reflect epidemiological reality \citep{Wearing2005,Greenhalgh2021}. Consequently, more general epidemic models have been developed to explicitly incorporate variable durations of infection and immunity, heterogeneous infectivity rates, or infection-age structure \citep{Metz1986,Inaba2001,Breda2012,Meehan2019}. These models increase biological realism and typically take the form of structured partial differential equations (PDEs) or renewal equations (REs) that offer two distinct, but equivalent, perspectives on the dynamics driving an epidemic.  

More precisely, structured PDEs track the dynamics of cohorts of infected individuals by explicitly modelling their progression from the initial infection to their recovery and potential return to the susceptible population. Within this framework, the structure variable typically represents the time-since-infection, or the \emph{infection-age}, of an infected individual, which is distinct from chronological time. In this cohort-centric perspective, the duration of infection can be modelled by a positive random variable which determines the recovery rate \citep{Cassidy2018a,Gurney1986,Feng2000}. It is also natural to assume that the infectivity of an infected individual changes as they progress through the infection, and thus include infection-age dependence in their contribution to the force of infection \citep{Kermack1927,Breda2012,Inaba2001,Greenhalgh2021}. On the other hand, REs model how an initial cohort of infected individuals propagates infection through successive generations \citep{Inaba2012a,Diekmann1990,Gyllenberg2007}, and this generational perspective reflects closely the original derivation by \citet{Kermack1927}. The RE, and thus the epidemic dynamics, are determined by combining the two biological processes that drive the infectious population: the production of new infected individuals and the duration of infection \citep{Breda2012,Diekmann1990}. 

Here, we leverage the equivalence between the PDE and the RE formulation to illustrate how including temporary immunity and variable infectiousness can induce nontrivial epidemic dynamics. We consider a minimal epidemic model that only tracks susceptible and infected individuals, and can therefore be considered an extension of the `standard' SIS model (that assumes constant infectivity and exponentially distributed infectious period).  We let both the infectivity and the rate of return to the susceptible class be functions of the time-since-infection. While our modelling approach does not explicitly include the recovered compartment typically used to model full immunity, it does allow for individuals to have a negligible contribution to the force of infection while remaining in the infected class. As these individuals are protected from re-infection, our model captures waning immunity through the rate at which individuals return to the susceptible class. In short, our modelling framework allows for the existence of a period in which individuals are not infectious but not yet susceptible. In this sense, we extend the SIS model considered by \citet{diekmann1982prelude} by incorporating waning immunity and infectiousness with a non-fixed duration.  

While our modelling framework accounts for general shapes of the infectivity and the probability distribution of the infectious/immune period, we consider gamma-distributed infectious periods with continuous shape parameter $j$. However, numerical methods are not commonly available for the resulting infinite-dimensional system (either in the PDE or RE formulation). We therefore compare three techniques to approximate the model with a finite-dimensional system of ODEs. We first consider an Erlang approximation, which consists of rounding the shape parameter $j$ to an integer value before using the linear chain trick \citep{Vogel1961,MacDonald1978,Cassidy2020a,Champredon2018,Diekmann2017,Ando2020}. We next consider a hypoexponential approximation, that is similar to the Erlang approximation, but preserves the mean and variance of the original gamma distribution~\citep{Cassidy2022}. Finally, we consider the pseudospectral approximation which is a discretization technique that works for more general distribution kernels~\cite{Breda2016,Gyllenberg2018,Scarabel2024Infinite}. By using these approaches to perform bifurcation analysis and time integration, we demonstrate how the inclusion of infection-age dependent parameters can result in dynamics that are not possible in the standard SIS model, including stable oscillations and genuine bistability between the endemic equilibrium and an oscillatory state. We then show how the common Erlang approximations can lead to a phenomenon of \emph{spurious bistability}, i.e., spurious dynamical outcomes that are purely determined by the chosen approximation technique, rather than the underlying model \citep{Humphries1993,Cassidy2022}.  

The genuine bistability results from the inclusion of more realistic sojourn distributions, as it is well known that the standard SIS and SIRS ODE models admit only a stable endemic equilibrium for $R_0 >1$, and no periodic oscillations. However, bistability is rarely observed in infection-age structured models as the analysis of periodic solutions beyond Hopf bifurcation points is technically challenging. Consequently, many existing models utilize the Erlang approximation by implicitly or explicitly assuming that the duration of infection is Erlang distributed, as this assumption allows for the model to be reduced to an equivalent system of ODEs. Indeed, \citet{hethcote1981nonlinear}, \citet{blyuss2010stability} and \citet{rost2020stability} showed that stable periodic solutions can emerge when waning of immunity is described by an Erlang distribution. Furthermore, \citet{taylor2009sir} showed that a delayed SIRS model with waning immunity and demographic turnover admits bistability near Hopf bifurcations, but, to our knowledge, this remains the only example in a model with waning immunity. More broadly, bistability is typically observed only in models that incorporate additional mechanisms, such as boosting of immunity upon re-exposure to the pathogen \citep{dafilis_frascoli_wood_mccaw_2012,kuhn,childs,opoku-sarkodie_dynamics_2022,opoku-sarkodie_bifurcation_2024,scarabel2025SIRS}. Here, we use numerical bifurcation analysis to show the existence of regions of bistability where a stable endemic equilibrium coexists with a stable periodic orbit in our minimal model that only includes waning immunity via the return to the susceptible class. Furthermore, the pseudospectral approximation enables us to study how the dynamical properties of the system change when varying the shape parameter, $j$, of the gamma-distributed infectious period in a continuous interval and shows the transition from stability to instability of the equilibrium.  This analysis not possible when using the linear chain trick and Erlang approximation, which enforces $j \in \mathbb{N}$. 
 
As mentioned, we explicitly consider infectious periods that are gamma, and not Erlang, distributed. While this distinction may appear trivial, it is relatively common to define an epidemiological model as a system of compartmental ODEs before fitting the model to disease incidence data \citep{Greenhalgh2021,Krylova2013}. The compartmental approach implicitly imposes an artificial relationship between the mean, $\tau$, and variance, $\sigma^2$, of the duration of infection, namely that $\tau^2= j \sigma^2$ for $j \in \mathbb{N}$. However, estimating the basic reproduction number, $R_0$, from outbreak data depends not only on the mean duration of infection, but also on the distribution of this period \citep{Roberts2007,Wearing2005,whittaker2023uncertainty}. Consequently, implicitly assuming that the duration of infection is Erlang, rather than gamma, distributed may lead to spurious estimates of $R_0$ and incorrectly predict the corresponding critical vaccination coverage. Accordingly, we here show how to relax the assumption $\tau^2= j \sigma^2$ while retaining the ability to numerically simulation and investigate an epidemic model.  

Overall, our results highlight the importance of distributional shape in determining epidemic dynamics. We show that the shape of the distribution of infection has a strong effect on the stability of the endemic equilibrium, even when the basic reproduction number $R_0$ and the mean duration of infectiousness are fixed. This implies that rounding or truncating shape parameters to integer values, as commonly done when using the linear chain trick, can lead to qualitatively different dynamics \citep{hethcote1981nonlinear,gonccalves2011oscillations,rost2020stability,Cassidy2022}. By combining RE theory, PDE formulation, and numerical bifurcation analysis, we provide new insights into the interplay between infectivity, waning immunity, and stability in infection-age structured models. Further, our analysis illustrates how the inclusion of biologically realistic mechanisms can lead to much richer dynamical behavior. 
 
\section{A simple SIS mathematical model structured by time-since-infection} 

We consider a minimal epidemic model for a constant population of $N$ individuals divided into two sub-populations representing individuals either susceptible to or infected by an infectious disease. As is commonly done, we consider the fraction of the total population who are susceptible, $S(t)$, or infected, $I(t)$. Consequently, we have $S(t) + I(t) = 1$. We begin by recalling the well-known SIS model without infection-age dependence. 

\subsection{An SIS model without infection-age dependence} 
In the well-known SIS ODE model, susceptible individuals are infected following contact with infected (and infectious) individuals with constant rate $\beta$. Infectious individuals recover with rate $\gamma$, and return to the susceptible compartment. The model is given by
\begin{equation}\label{Eq:SISModelODE}
\left.
    \begin{aligned}
        \TimeDeriv S(t)   = \gamma I -\beta S I \\
        \TimeDeriv I(t) = \beta S I - \gamma I. 
    \end{aligned}
    \right \}
\end{equation}
In this model, the duration of infection is exponentially distributed with mean duration $1/\gamma$, and individuals are infectious with constant infectivity rate $\beta$ throughout infection~\cite[Section 2.3]{keeling2008modeling}.  

We now extend this simple model framework to allow for the infectivity, $\beta$, and the duration of infection to depend on the time passed since infection. We begin from the perspective of tracking cohorts of individuals to derive a structured PDE model before developing the equivalent RE model from the generational perspective. 

\subsection{An SIS model with infection-age dependence: the cohort perspective} 

We denote the density of infected individuals with infection-age $a$ at time $t$, i.e. those who were infected at time $t-a$ and have not returned to the susceptible class, by $i(t,a).$ We calculate the total fraction of infected individuals by integrating over all possible infection-ages
\begin{align}\label{Eq:totalInfected}
I(t) = \int_0^{\infty} i(t,a) \d a. 
\end{align}
Now, recalling that the total population is closed, so that $S(t) + I(t) = 1$, the total susceptible population is given by 
\begin{align*}
S(t) = 1 - \int_0^{\infty} i(t,a) \d a = 1-I(t).
\end{align*}
We assume that all individuals are equally susceptible to infection and therefore do not include any structure in the susceptible population, although it is possible to include differing levels of susceptibility to infection \citep{Inaba2001,Meehan2020}. Then, assuming that the population is well-mixed so that infection follows mass action kinetics, new infections occur proportionally to the number of susceptible individuals via 
\begin{align*}
i(t,0) = \Lambda(t) S(t).
\end{align*} 
Here, $\Lambda(t)$ is the force of infection, i.e., the per capita probability per unit time of a susceptible individual to acquire infection. We assume that an infected individual with infection-age $a$ has infectiousness $\beta(a)$, so that $\Lambda(t)$ is given by  
\begin{equation}
\Lambda(t) = \int_0^{\infty} \beta(a)i(t,a) \d a.
\label{Eq:ForceOfInfection}
\end{equation}
We impose that $\beta(a) \geq 0 $ and that $\beta \in L_1(\mathbb{R}_+)$. Note that $\beta$ can have compact support if the infectious period is uniformly bounded. 

Now, to close the model, we must describe how infected individuals return to the susceptible compartment. We assume that the duration of infection is a positive random variable $\mathbb{A}$ with probability density function (PDF) $K$, such that 
\begin{align*}
K(a) \geq 0\quad \textrm{and} \quad \int_0^{\infty} K(a) \d a = 1. 
\end{align*}
The rate at which infected individuals return to the susceptible class is given by the hazard rate \citep{Cassidy2018a,Gurney1986} 
\begin{equation}
h(a) = \lim \limits_{\delta a \to 0} \frac{P \left[ a< \mathbb{A} <a+\delta a \right] }{\delta a} = \frac{K(a)}{1-\int_0^a K(s) \d s}. 
\label{Eq:HazardRate}
\end{equation}
As mentioned, the rate at which these individuals return to the susceptible class corresponds to the loss of protection from re-infection and thus waning immunity. We can therefore describe the dynamics of the infected class by the integro-partial differential equation 
\begin{equation} \label{Eq:InfectedStructuredPDEandBC} 
\left.
\begin{aligned}
& \left( \partial_t  + \partial_a \right) i(t,a) = -h(a) i(t,a)  \\
& i(t,0) = \Lambda(t) S(t) = \int_0^\infty \beta(a)i(t,a) \d a \left( 1 - \int_0^\infty i(t,a) \d a \right) .
\end{aligned} 
\right \}
\end{equation}
The left-hand side of the PDE~\eqref{Eq:InfectedStructuredPDEandBC} describes the progression of an individual through the infected compartment, while the right-hand side models the rate at which these individuals, with infection-age $a$, return to the susceptible class. The boundary condition of Eq.~\eqref{Eq:InfectedStructuredPDEandBC}, $i(t,0),$ describes the inflow of newly infected individuals due to new infections at time $t$. 

\begin{figure}[tp]
    \centering
\includegraphics[trim= 4 10 5 10,clip,width=1\textwidth]{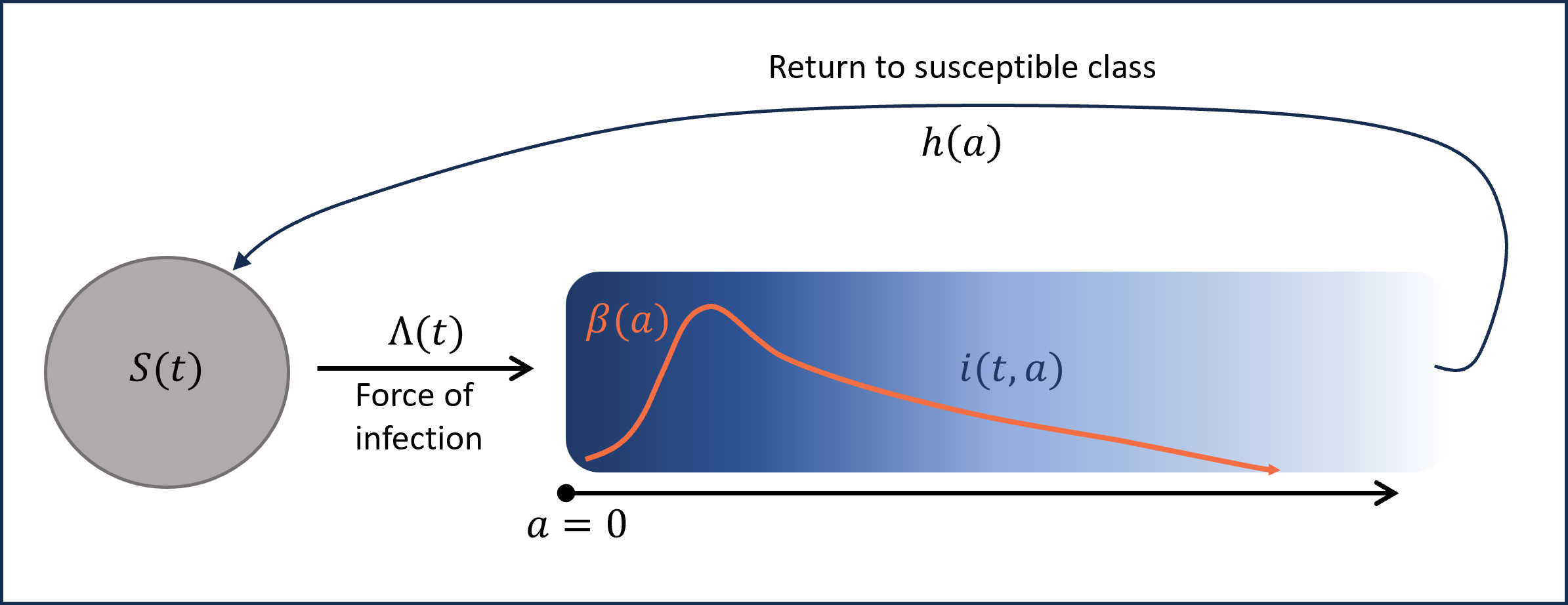}
    \caption{\textbf{ A model schematic of the SIS model structured by infection-age}. Susceptible individuals are denoted by $S(t)$ and are infected with force of infection $\Lambda(t)$. Newly infected individuals have infection-age $a=0$ and enter at the left-hand boundary of the age-structured infected population. As these infected individuals progress through infection, their infectivity, and corresponding contribution to the force of infection, evolves according to $\beta(a) $. Infected individuals return to the susceptible class at infection-age dependent rate $h(a),$ which is the hazard rate corresponding to the random variable $\mathbb{A}$ that defines the duration of infection and corresponds to the immunity waning rate. Individuals can remain in the infectious class but no longer contribute to the force of infection if $\beta(a)$ has compact support. }
    \label{Fig:ModelDiagram}
\end{figure}

As mentioned, we decouple the infectivity rate, $\beta(a)$ from the rate at which individuals return to the susceptible class, $h(a)$. Consequently, this simple SIS model can capture a `recovered' phase, where individuals are no longer infectious with $\beta(a) = 0,$ but have not yet returned to the susceptible class, and are thus immune to re-infection. We therefore do not include an explicit `recovered' compartment in our model to reduce the number of nested integrals in the model formulation. We give a schematic of the model formulation in Figure~\ref{Fig:ModelDiagram}. 

\subsection{An SIS model with infection-age dependence: the generational perspective} 
We now derive the equivalent RE formulation of Eq.~\eqref{Eq:InfectedStructuredPDEandBC}. Along the characteristics of Eq.~\eqref{Eq:InfectedStructuredPDEandBC}, the density of infected individuals is given by $i(t,a) = i(t-a,0) \mathcal{F}(a)$, where 
\begin{align}
\mathcal{F}(a) = 1- \int_{0}^a K(s) \d s 
\end{align}
is the survival probability in the infected compartment, i.e., the probability that an individual has not exited the infected compartment by infection-age $a$. Substituting this expression for $i(t,a)$ into the boundary condition of Eq.~\eqref{Eq:InfectedStructuredPDEandBC} gives the equivalent nonlinear RE for $i(t,0)$
\begin{align}
i(t,0) & = \Lambda(t) S(t) = \left(1-\int_0^\infty i(t-a,0)\mathcal{F}(a) \d a \right) \int_0^{\infty} \beta(a) \mathcal{F}(a) i(t-a,0) \d a. 
\label{Eq:InfectedRenewalEquation}
\end{align}

However, it is possible to derive Eq.~\eqref{Eq:InfectedRenewalEquation} directly from the modelling ingredients by noting that the kernel $A(a)=\beta(a) \mathcal{F}(a)$ describes the effective contribution to the force of infection of an individual of infection-age $a$ \citep{Diekmann1990,Inaba2012,Inaba2008,DiekmannScarabelAge}, and $\mathcal{F}(a)$ describes the probability that an individual has not returned to the susceptible class at infection-age $a$. Here, we have presented both formulations to explicitly highlight the roles of the two processes included in the model, infection and return to the susceptible class, and link the generational and cohort perspectives. 

  
A solution to Eq.~\eqref{Eq:InfectedStructuredPDEandBC} for $t\geq 0$ can be uniquely determined by specifying an initial condition $i(0,a)=\phi_0(a)$ for $a \in (0,\infty)$, with $\phi_0\geq 0$ an $L^1$ function satisfying appropriate conditions. This can be uniquely mapped to an initial condition for Eq.~\eqref{Eq:InfectedRenewalEquation} for $t<0$ (the initial `history') by setting 
\begin{align*}
\phi_0(a) = i(-a,0) \mathcal{F}(a), \quad a>0.
\end{align*}
Effectively, this definition of $i(-a,0)$ traces the initial cohort of infected individuals back to their infection time $-a<0$~\citep{Cassidy2018a,DiekmannScarabelAge}.

\subsection{Equilibria, stability, the Malthusian parameter and the basic reproduction number}

We now focus on the long-term dynamics of the system and assume that model solutions are globally defined for $t \in \mathbb{R}$. We immediately see that Eq.~\eqref{Eq:InfectedRenewalEquation} admits the disease-free equilibrium (DFE), where $S\equiv 1$ and $i(t,0) \equiv 0$.  The principle of linearised stability for nonlinear REs implies that the local stability of this equilibrium is determined by the roots of the characteristic equation obtained from the linearisation at the corresponding equilibrium \citep{Diekmann2012}. The linearisation of Eq.~\eqref{Eq:InfectedRenewalEquation} at an equilibrium $i(t,0)\equiv x^*$ reads
\begin{equation}
\label{Eq:LinearisedGeneral}
    x(t) = - x^* R_0 \int_0^\infty x(t-a)\mathcal{F}(a)\d a + \left( 1-x^*\int_0^\infty \mathcal{F}(a)\d a\right) \int_0^\infty x(t-a)\beta(a)\mathcal{F}(a)\d a. 
\end{equation}
Now, at the DFE, Eq.~\eqref{Eq:LinearisedGeneral} reduces to 
\begin{equation*}
    x(t) = \int_0^\infty x(t-a)\beta(a)\mathcal{F}(a)\d a,
\end{equation*}
and, using the trial solution $x(t) = e^{\lambda t}$ for $\lambda \in \mathbb{C}$, the corresponding characteristic equation is 
\begin{equation} \label{Eq:EulerLotka}
    1 = \int_0^\infty e^{-\lambda a} \beta(a)\mathcal{F}(a)\d a =: F(\lambda). 
\end{equation}
Eq.~\eqref{Eq:EulerLotka} is typically known as the Euler--Lotka equation \cite{Diekmann2012a}. Since the right-hand side $F(\lambda)$ is a decreasing function of $\lambda$, Eq.~\eqref{Eq:EulerLotka} can have only one real solution $\lambda=r$ and all complex roots of Eq.~\eqref{Eq:EulerLotka} have smaller real part than $r$~\citep{Metz1986,Diekmann1990}. The unique real root $r$ is the intrinsic growth rate of the infected population at the beginning of an epidemic, and is known as the Malthusian parameter. We note that the Malthusian parameter can be obtained from the PDE in Eq.~\eqref{Eq:InfectedStructuredPDEandBC} \citep{Perthame2006}. The sign of the Malthusian parameter, and thus the local stability of the DFE, is determined by the basic reproduction number, 
\begin{equation}\label{Eq:R0}
R_0 = F(0) = \int_0^\infty \beta(a) \mathcal{F}(a)\d a, 
\end{equation}
that describes the expected number of secondary infections generated by one infected individuals throughout their infectious period, in a completely susceptible population \cite{Diekmann1990}. Consequently, the classic sign relationship $\textrm{sign}\left( R_0-1 \right) = \textrm{sign}(r)$ holds: $r$ is negative and DFE is stable if $R_0<1$, while $r$ is positive and DFE is unstable if $R_0>1$. 

Now, if $R_0>1$, there is a positive endemic equilibrium (EE) given by  
$$ (EE) \quad S\equiv \frac{1}{R_0} \quad \text{and}\quad i(t,0)\equiv \left(1-\frac{1}{R_0}\right) \frac{1}{\int_0^\infty \mathcal{F}(a)\d a}. $$
The characteristic equation corresponding to the EE is 
\begin{equation*}
    1 = \frac{\left( 1-R_0\right)}{\int_0^\infty \mathcal{F}(s)\d s} \int_0^\infty e^{-\lambda a} \mathcal{F}(a) \d a + \frac{1}{R_0} \int_0^\infty e^{-\lambda a} \beta(a)\mathcal{F}(a) \d a. 
\end{equation*}
Consequently, the stability of the EE depends on the Laplace transforms of the functions $\mathcal{F}$ and $\beta\mathcal{F}$. The analysis of the roots of this equation is significantly more complex than Eq.~\eqref{Eq:EulerLotka}. We will show via numerical investigations that EE can undergo a Hopf bifurcation leading to stable periodic solutions, that may coexist with the stable EE. In particular, we will show that the existence and stability of periodic solutions depend not only on $R_0$, but also on the shape of the functions $\beta$ and $\mathcal{F}$.

\section{Ordinary differential equation approximations of the epidemic dynamics} \label{Sec:ODEApproximations}

In this section, we derive ODE approximations of the simple epidemic model. We start by considering the total number of infected individuals, given by
\begin{equation*}
I(t) = \int_0^{\infty} i(t,a) \d a = \int_{-\infty}^t i(t,t-s) \d s.
\end{equation*}
Differentiating $I(t)$ using Leibniz's rule gives
\begin{align}\label{Eq:InfectedDistributeDDE}
\TimeDeriv I(t) &  = \Lambda(t)S(t) - \int_{0}^{\infty}   \Lambda(t-a)S(t-a) K(a) \d a, 
\end{align}
where we have used the PDE Eq.~\eqref{Eq:InfectedStructuredPDEandBC}, the solution $i(t,a)=i(t-a,0)\mathcal{F}(a)$, and Eq.~\eqref{Eq:HazardRate}. Now, while we can write $S(t) = 1-I(t)$, the differential equation for $I(t)$ also depends on the force of infection $\Lambda(t)$. We therefore use a similar technique to obtain a differential equation for $\Lambda(t)$. Once again using Leibniz's rule (and assuming the derivative $\beta'$ exists), we calculate
\begin{align} \label{Eq:LambdaGenericDE-GenericBeta}
\TimeDeriv \Lambda(t) & =  \beta(0)\Lambda(t)S(t) + \int_0^\infty \Lambda(t-a)S(t-a) \left[ \beta'(a) \mathcal{F}(a) - \beta(a) K(a) \right] \d a. 
\end{align} 
We thus obtained two distributed delay differential equations (DDEs) for $I(t)$ and $\Lambda(t)$, where the delayed term is given by the convolution integrals. 

The DDEs Eq.~\eqref{Eq:InfectedDistributeDDE} and Eq.~\eqref{Eq:LambdaGenericDE-GenericBeta} should be paired with appropriate initial conditions for $t\leq 0$. Given an initial condition $i(0,a)=\phi_0(a)$ for Eq.~\eqref{Eq:InfectedStructuredPDEandBC}, we define the history functions for $I$ and $\Lambda$ by 
\begin{align}
I(s) &= \int_0^\infty i(s,a) \d a = \int_0^{\infty} i(s-a,0)\mathcal{F}(a)\d a  = \int_0^{\infty} \phi_0(a-s) \frac{\mathcal{F}(a)}{\mathcal{F}(a-s)}\d a, \label{Eq:IHistory}
\end{align}
and 
\begin{align}
\Lambda(s) &= \int_0^{\infty} i(s-a,0)\beta(a) \mathcal{F}(a)\d a  = \int_0^{\infty} \phi_0(a-s) \beta(a) \frac{\mathcal{F}(a)}{\mathcal{F}(a-s)}\d a, \label{Eq:LambdaHistory}
\end{align}
for $s\leq 0$. 

We have thus far not chosen specific forms for $K(a)$ or $\beta(a)$. To simplify the analysis that follows, we take 
\begin{align}\label{Eq:ForceInfectionDefinition}
\beta(a) = \beta_0 e^{-k_d a}, \quad a\geq 0, 
\end{align}
where $\beta_0>0$ is the maximal infectivity of an infectious individual and $k_d>0$ is the rate at which this infectivity wanes. However, we note that other choices may be appropriate when attempting to fit Eq.~\eqref{Eq:InfectedRenewalEquation} to data.
With this choice and using Eq.~\eqref{Eq:ForceOfInfection}, the DDE Eq.~\eqref{Eq:LambdaGenericDE-GenericBeta} further reduces to 
 \begin{align} \label{Eq:LambdaGenericDE}
\TimeDeriv \Lambda(t)  & =  \beta_0 \Lambda(t)S(t) - k_d\Lambda(t)  - \int_{0}^{\infty}   \Lambda(t-a)S(t-a) \beta(a) K(a)  \d a.
\end{align} 

 Now, if $K$ is an Erlang distribution (i.e., a gamma distribution with integer shape parameter), the DDEs \eqref{Eq:InfectedDistributeDDE} and \eqref{Eq:LambdaGenericDE} can be reformulated as mathematically equivalent systems of ODEs via the linear chain trick. However, the PDF $K$ can be estimated from epidemic data, and an exact ODE reduction is, in general, impossible for other kernels  \citep{Diekmann2017,Diekmann2020}. In what follows, we assume that the duration of infection is gamma distributed with mean $\tau$ and variance $\sigma^2$. The PDF of the gamma distribution is given by
\begin{align*}
g_b^j(s) = \frac{b^{j} s^{j-1}e^{-bs}}{\Gamma(j-1)}, \quad s\geq 0, 
\end{align*}
where the mean and variance determine the rate, $b$, and shape, $j$, parameters of the PDF $K$ via
 \begin{align*}
     \tau = \frac{j}{b} \quad \textrm{and} \quad \sigma^2 = j/b^2. 
 \end{align*}
As mentioned, we allow the shape parameter $j$ to be any positive real, rather than imposing $j \in \mathbb{N}$ \citep{Champredon2018,Inaba2008,Hurtado2025}. 

In the rest of this section, we present three mathematical approaches to reduce the epidemic model to a system of ODEs. We first follow \citet{Cassidy2020a} and derive a system of ODEs whose solution approximates the dynamics of the epidemic modelled by Eqs.~\eqref{Eq:InfectedDistributeDDE} and \eqref{Eq:LambdaGenericDE} through an application of the linear chain trick \citep{deSouza2017,Cassidy2018,Cassidy2018a}. We next present an improved approximation via hypoexponential distributions \citep{Cassidy2022}. It is important to note that, if the gamma distribution is indeed an Erlang distribution, both of these approximations are exact and the corresponding ODE system is mathematically equivalent to Eqs.~\eqref{Eq:InfectedDistributeDDE} and \eqref{Eq:LambdaGenericDE}. In these approximations, the dimension of the approximating ODE system is determined by the kernel $K$. Finally, we summarize a numerical approximation of the RE \eqref{Eq:InfectedRenewalEquation}, which applies to general kernels $K$ and reduces the RE to a system of ODEs via pseudospectral discretization. In the pseudospectral case, the dimension of the approximating ODE system is a numerical, rather than model, parameter. 

\subsection{Erlang approximation}

The Erlang approximation of a generic gamma distribution with mean $\tau $ and variance $\sigma^2$ is well known within the mathematical community as the \textit{linear chain trick} \citep{Vogel1961,MacDonald1978}. Specifically, this approach allows modellers to replace the convolution integral in Eq.~\eqref{Eq:InfectedDistributeDDE} by the solution of a system of ODEs. However, the requirement that $\tau^2$ is an integer multiple of $\sigma^2$ is not generic, so modellers typically round $j$ to the nearest integer $[j]$ (with $[0.5]=1$). Then, setting $b = [j]/\tau$ ensures that the Erlang approximation has the same mean, $\tau$, as the underlying gamma distribution. However, due to the rounding of $j$, the variance of the resulting Erlang distribution is $(j^2/[j]^2)\sigma^2$. 

After approximating the gamma distribution by an Erlang distribution, the \textit{Erlang approximation} of Eq.~\eqref{Eq:InfectedDistributeDDE} is
\begin{align}\label{Eq:ErlangApproximationInfectedDistributeDDE} 
\TimeDeriv I(t) &  = \Lambda(t)S(t) - \int_{0}^{\infty}   \Lambda(t-a)S(t-a) g^{[j]}_b(a) \d a. 
\end{align}
Now, we define
\begin{align} \label{Eq:InE-Erlang}
    I_{n,E}(t) = \int_{0}^{\infty}   \Lambda(t-a)S(t-a) g^{n}_b(a) \d a, \quad \textrm{for} \quad n = 1, 2,3,\dots,[j],
\end{align}
so that the convolution integral in Eq.~\eqref{Eq:ErlangApproximationInfectedDistributeDDE} is precisely $I_{[j],E}(t)$, where the subscript $E$ denotes the Erlang approximation. Now, the PDF of the Erlang distribution satisfies
\begin{align}\label{Eq:ErlangDerivative}
    \frac{\d }{\d s} g_b^{1}(s) = -bg_b^1(s) , \quad \textrm{and} \quad  \frac{\d }{\d s} g_b^{i}(s) = b\left( g_b^{i-1}(s)- g_b^{i}(s) \right) \quad \textrm{if} \quad i = 2,3,\dots.
\end{align}
Consequently, differentiating $I_{[j],E}(t)$ using Leibniz's rule and using Eq.~\eqref{Eq:ErlangDerivative}, we obtain
\begin{equation}\label{Eq:InfectedErlangDE}
\left. 
\begin{aligned}
    \TimeDeriv I_{1,E}(t) & = \Lambda(t)S(t) - bI_{1,E}(t) \\
    \TimeDeriv I_{i,E}(t) & = b(I_{i-1,E}(t) - I_{i,E}(t) ) \quad \textrm{for} \quad i = 2,3,\dots,[j].
\end{aligned}
\right \}
\end{equation}
We next turn to Eq.~\eqref{Eq:LambdaGenericDE} and define 
\begin{align} \label{Eq:LnE-Erlang}
    L_{n,E}(t) = \int_{0}^{\infty} \Lambda(t-a)S(t-a) \beta(a)  g^{n}_b(a) \quad \textrm{for} \quad n = 1,2,3,\dots,[j].
\end{align}
Then, a similar approach yields \citep{Cassidy2020a,Cassidy2018} 
\begin{equation} \label{Eq:LambdaErlangDE}
\left. 
\begin{aligned}
\TimeDeriv L_{1,E}(t) & = \beta_0 \Lambda(t)S(t) - \left(b+k_d \right)  L_{1,E}(t) \\
\TimeDeriv L_{n,E}(t) & =  b \left[ L_{n-1,E}(t) - L_{n,E}(t) \right] - k_d L_{n,E}(t), \quad \textrm{for} \quad n = 2,3,\dots,[j]. 
\end{aligned} 
\right \}
\end{equation}
Consequently, we can approximate the convolution integrals in the gamma distributed DDEs Eq.~\eqref{Eq:InfectedDistributeDDE} and \eqref{Eq:LambdaGenericDE} by the systems of ODEs in Eq.~\eqref{Eq:InfectedErlangDE} and Eq.~\eqref{Eq:LambdaErlangDE}, respectively. We note that this approach can be extended to the case when the function $\beta$ is itself a scaled Erlang distribution.

\subsection{Hypoexponential approximation}

As mentioned, the Erlang approximation involves replacing the gamma distribution by an Erlang distribution that precisely matches the mean, or first moment, of the underlying gamma distribution. \citet{Cassidy2022} showed that it is possible to develop a hypoexponential distribution that precisely matches the first and second moment of the underlying gamma distribution. Specifically, they proved the following result (where $\ceil{\cdot}$ denotes the ceiling operator). 
\begin{theorem}\label{Thm:TwoRateApprox}[Theorem~3.1 of \citet{Cassidy2022}]
Consider the random variable $\mathcal{X}$ with mean $\tau$ and variance $\sigma^2$, with shape parameter $j = \tau^2/\sigma^2$. Let $\mathcal{Y}$ be the random variable obtained by concatenating $n = \max( \ceil{j},2 )$ independent and exponentially distributed random variables where $n-2$ of these exponentially distributed variables have identical rates
\begin{align*}
\lambda_{i} = \lambda = \frac{n}{\tau}\,,\quad i = 1,\dots, n-2\,,
\end{align*}
while the remaining two exponentially distributed variables have rates $\lambda_{n-1} = \nu$ and $ \lambda_{n} = \mu$.  Then, setting 
\begin{align*}
\nu = \left( \frac{\tau}{n}\left(1+\sqrt{\tfrac{n}{2j}(n-j)}\right)  \right)^{-1}  
\end{align*} 
and
\begin{align*}
\mu = \left(  \frac{\tau}{n}\left(1-\sqrt{\tfrac{n}{2j}(n-j)}\right) \right)^{-1}
\end{align*}
ensures that $\mathcal{X}$ and $\mathcal{Y}$ have the same first two moments.
\end{theorem}
The random variable $\mathcal{Y}$ is the concatenation of $n$ exponentially distributed random variables and thus a hypoexponential distribution. Now, we use this hypoexponential distribution to approximate the convolution integrals in Eq.~\eqref{Eq:InfectedDistributeDDE} and Eq.~\eqref{Eq:LambdaGenericDE}. Here, the hypoexponential approximation to the convolution integral is given by 
\begin{align} \label{Eq:InE-hypo}
    I_{H}(t) = \int_{0}^{\infty}   \Lambda(t-a)S(t-a) f_{\mathcal{Y}}(a) \d a.
\end{align}
Now, recalling that $\mathcal{Y}$ is a concatenation of $n$ exponential distributions, we define the auxiliary variables $I_{m,H}(t)$ that satisfy 
\begin{equation} \label{Eq:InfectedHypoODE}
    \left. 
\begin{aligned}
    & \TimeDeriv I_{1,H}(t)  = \Lambda(t)S(t) - bI_{1,H}(t) \\
    & \TimeDeriv I_{i,H}(t)  = \lambda (I_{i-1,H}(t) - I_{i,H}(t) ) \quad \textrm{for} \quad i = 2,3,\dots,n-2 \\
    & \TimeDeriv I_{n-1,H}(t)  = \lambda I_{n-2,H}(t) - \nu I_{n-1,H}(t) \\
    & \TimeDeriv I_{H}(t)  =  \nu I_{n-1,H}(t) - \mu I_{H}(t).
\end{aligned}
\right \}
\end{equation}
Similarly, we define the hypoexponential approximation to the convolution integral in Eq.~\eqref{Eq:LambdaGenericDE} by
\begin{align} \label{Eq:LnE-hypo}
    L_{H}(t) = \int_{0}^{\infty} \Lambda(t-a)S(t-a) \beta(a) f_{\mathcal{Y}}(a) \d a  ,
\end{align}
and we obtain 
\begin{equation}\label{Eq:LambdaHypoODE}
\left.
\begin{aligned}
& \TimeDeriv L_{1,H}(t)   = \beta_0 \Lambda(t)S(t) - \left(\lambda+k_d \right)  L_{1,H}(t) \\
& \TimeDeriv L_{i,H}(t)   =  \lambda \left[ L_{i-1,H}(t) - L_{i,H}(t) \right] - k_d L_{i,H}(t), \quad \textrm{for} \quad i = 2,3,\dots, n-2 \\
& \TimeDeriv L_{n-1,H}(t) = \lambda L_{n-2,H}(t) - \nu L_{n-1,H}(t) - k_d L_{n-1,H}(t) \\
& \TimeDeriv L_{H}(t) = \nu L_{n-1,H}(t) - \mu L_{H}(t) - k_d L_{H}(t) .
\end{aligned}
\right \}
\end{equation}
Thus, we can once againapproximate the convolution integrals in the gamma distributed DDEs Eq.~\eqref{Eq:InfectedDistributeDDE} and \eqref{Eq:LambdaGenericDE} by the systems of ODEs in Eq.~\eqref{Eq:InfectedHypoODE} and Eq.~\eqref{Eq:LambdaHypoODE}, respectively, and this approach can be generalised to the case when the function $\beta$ is a scaled Erlang kernel. 

Now, the initial conditions for both the Erlang and hypoexponential approximations are immediately determined by the initial distribution $\phi_0$. The initial fraction of infected individuals and force of infections are given by  
\begin{align*}
I(0) = \int_0^{\infty} \phi_0(a) \d a,  \quad \textrm{and} \quad \Lambda(0) = \int_0^{\infty} \phi_0(a) \beta(a) \d a,
\end{align*}
while the initial conditions for the auxiliary variables, $I_{n,\cdot}(0)$ and $L_{n,\cdot}(0)$ are then prescribed via the history of $I$ and $\Lambda$ using Eqs.~\eqref{Eq:InE-Erlang} and \eqref{Eq:LnE-Erlang} (for the Erlang approximation) or Eqs.~\eqref{Eq:InE-hypo} and \eqref{Eq:LnE-hypo} (for the hypoexponential approximation) \citep{Cassidy2018a,Cassidy2018,Cassidy2025}.

\subsection{The pseudospectral approximation for bifurcation analysis with general kernels}
A recently proposed numerical method derives an approximating system of ODEs via \emph{pseudospectral}, or \emph{spectral collocation}, methods, for general parameters $\beta$ and $K$. We here summarise the core ideas of this approach, while we refer to \citet{Gyllenberg2018, Scarabel2024Infinite,Scarabel2025Truncated} and references therein for further details and rigorous error bounds, and to \citet[Appendix A]{DiekmannScarabelSize} for a concise and less technical summary. 

Let $u(t)$ be the history of the solution $i(t,0)$ of Eq.~\eqref{Eq:InfectedRenewalEquation} at time $t$, namely 
\begin{equation*}
    u(t)(s) = i(t+s,0), \quad s \leq 0. 
\end{equation*}
Following \citet{Diekmann2012}, we assume there exists a parameter $\rho>0$ such that, for each $t \geq 0$, $u(t)$ belongs to the space of exponentially weighted integrable functions 
\begin{equation*}
    L^1_\rho((-\infty,0],\mathbb{R}) := \left\{ f\colon (-\infty,0] \to \mathbb{R} \text{ such that } \int_{-\infty}^0 |f(s)| e^{\rho s} \d s < \infty \right\}. 
\end{equation*}

To obtain a numerical approximation of the state $u(t)$, it is convenient to map it to a more regular function space via integration, and consider instead the state 
\begin{equation*}
    v(t)(\theta) = e^{\rho \theta} \int_0^{\theta} u(t)(s) \d s, \quad \theta \leq 0, 
\end{equation*}
that is absolutely continuous and tends to zero when $\theta\to -\infty$. 

Eq.~\eqref{Eq:InfectedRenewalEquation} can be reformulated as a semilinear abstract differential equation for $v(t) \in AC(\mathbb{R}_{\leq 0},\mathbb{R})$ as  
\begin{equation}\label{Eq:ADE}
v'(t)= A v(t) + F(v(t)), \quad t\geq 0, 
\end{equation}
where $A$ is a linear operator defined, using the weight $w(\theta) = e^{\rho\theta}$, by 
\begin{align*}
A(w\phi) &= w\phi', \quad w\phi \in D(A), \\
D(A) &= \{w\phi\ \mid\ w\phi,\ w\phi'\in AC(\mathbb{R}_{\leq 0},\mathbb{R}),\ \phi(0)=0 \},
\end{align*}
and $F$ is a nonlinear operator that captures the action of the right-hand side of the RE, given by 
\begin{equation*}\label{F}
F(w\phi) = - w(\cdot) \times \left( 1 - \int_0^\infty \phi'(-a)\mathcal{F}(a) \d a\right) \int_0^\infty \beta(a)\mathcal{F}(a) \phi'(-a) \d a . 
\end{equation*}

The main idea of spectral collocation is to approximate, for each $t$, the function $v(t)$ using an exponentially weighted polynomial on $\mathbb{R}_{\leq 0}$ with a given degree $d \in \mathbb{N}$ . By requiring that the weighted polynomials satisfy the abstract differential equation \eqref{Eq:ADE} on a given set of collocation nodes 
$\{\theta_1,\dots \theta_d\}$ with 
$$\theta_d<\theta_{d-1}<\cdots< \theta_1 < 0,$$ 
the RE is approximated by a system of $d$ ODEs~\citep{Scarabel2024Infinite}, where each variable represents an approximation of $v(t)$ at each collocation node. 

The resulting approximating system of ODEs takes the form 
\begin{equation*}
\TimeDeriv V(t)=A_d V(t)+F_d(V(t)), 
\end{equation*}
for $V(t) \in \mathbb{R}^d$, where 
$$V_n(t) \approx v(t)(\theta_n) = e^{\rho\theta_n} \int_0^{\theta_n} i(t+s,0)\d s, \quad n=1\dots,d.$$ 
The matrix $A_d \in \mathbb{R}^{d\times d}$ is the weighted differentiation matrix and is completely defined by the nodes and the parameter $\rho$. The nonlinear function $F_d$ captures the action of $F$: it applies $F$ to the weighted $(d+1)$-degree polynomial interpolating $V(t)$ on $\theta_1,\dots,\theta_d$, and zero in $\theta_0=0$. 
When integrals need to be approximated numerically, $F_d$ also accounts for the use of quadrature rules.
We obtain an approximation of $i(t,0)$ by applying the function $F_d$ to $V(t)$, so that:
$$ i(t,0) \approx F_d(V(t)). $$

A crucial aspect for an accurate approximation of the stability properties of equilibria is the choice of the collocation nodes in $\mathbb{R}_{\leq 0}$. We follow \citet{Scarabel2024Infinite} and take the extrema of the classical Laguerre orthogonal polynomials scaled by $1/(2\rho)$ \citep{Mastroianni}. The corresponding Gauss--Radau quadrature rules are then used to approximate the integrals on $\mathbb{R}_{\leq 0}$.  The numerical approximation can be computed with numerically stable algorithms for weighted polynomial interpolation and weighted pseudospectral differentiation~\citep{Mastroianni, weideman2000}. 
A simpler (in the sense of requiring less parameters) approximation based on Chebyshev nodes can be used when the kernels have compact support~\cite{Breda2016,Scarabel2021}. 


An important difference with the approximation approaches presented in the previous sections is that in this case the dimension $d$ of the ODE system is a numerical choice and is not determined by the model parameters. As $d$ increases, \citet{Scarabel2024Infinite} showed that the approximating system of ODEs reproduces exactly the equilibria of the RE \eqref{Eq:InfectedRenewalEquation}, and the approximate roots converge to the characteristic roots exponentially with respect to $d$. The ODE approximation is therefore reliable to study the stability and bifurcations of equilibria. Moreover, the ODE system can be conveniently studied with widely available software for ODEs. 
We note that the convergence of the solutions of the initial value problem is an open problem, so in the following we will use the pseudospectral approximation only to study the long-term properties of the systems. 

\section{Numerical bifurcation analysis (and time integration)} 

We numerically investigate  the stability of the EE using the pseudospectral approximation approach. We consider the infection-age dependent infectivity rate $\beta(a)$ given by Eq.~\eqref{Eq:ForceInfectionDefinition} and assume that the duration of infection is gamma-distributed with shape parameter $j>0$ and mean $\tau$, where we explicitly do not assume that $j\in \mathbb{N}$. Now, the basic reproduction number $R_0$, given in Eq.~\eqref{Eq:R0}, depends on both $\beta$ and $K$. However, $R_0$ is often used as a proxy for the infectivity of a pathogen and is of primary importance in public health decision making. For this reason, we consider $R_0$ as a parameter and compute the corresponding $\beta_0$ from the other parameters in our analyses. 

In what follows, we consider the mean duration of infection, $\tau$, and the rate at which infectivity decays motivated by respiratory syncytial virus (RSV). 
We take the mean duration of protection from infection of $\tau = 5$ months, which is consistent with the estimates in a recent systematic review \citep{Lang2022}. Further, in a prospective study of 47 households with 493 individuals, \citet{Munywoki2015} found the average duration of virus shedding following RSV infection was 11.2 days, which is similar to rates typically used for the duration of infectivity \citep{Lang2022}. We therefore take $k_d = 2,$ which implies that infectivity will have dropped by roughly $63\%$ in the two weeks post-infection. 

We perform numerical bifurcation analyses of the approximating ODE systems using MatCont (version 7p6), an established MATLAB software package for numerical continuation and bifurcation analysis that has been used throughout mathematical biology applications \citep{MatCont2008, Liessi-matcont, Cassidy2025}. The numerical bifurcation diagrams are obtained with $d=30$. The parameter $\rho$, used to define the collocation nodes in the pseudospectral approximation, is $2j/\tau$. 

\subsection{Bifurcation structure}

\begin{figure}[htp]
    \centering
    \includegraphics[trim= 2 10 5 15,clip,width=1\textwidth]{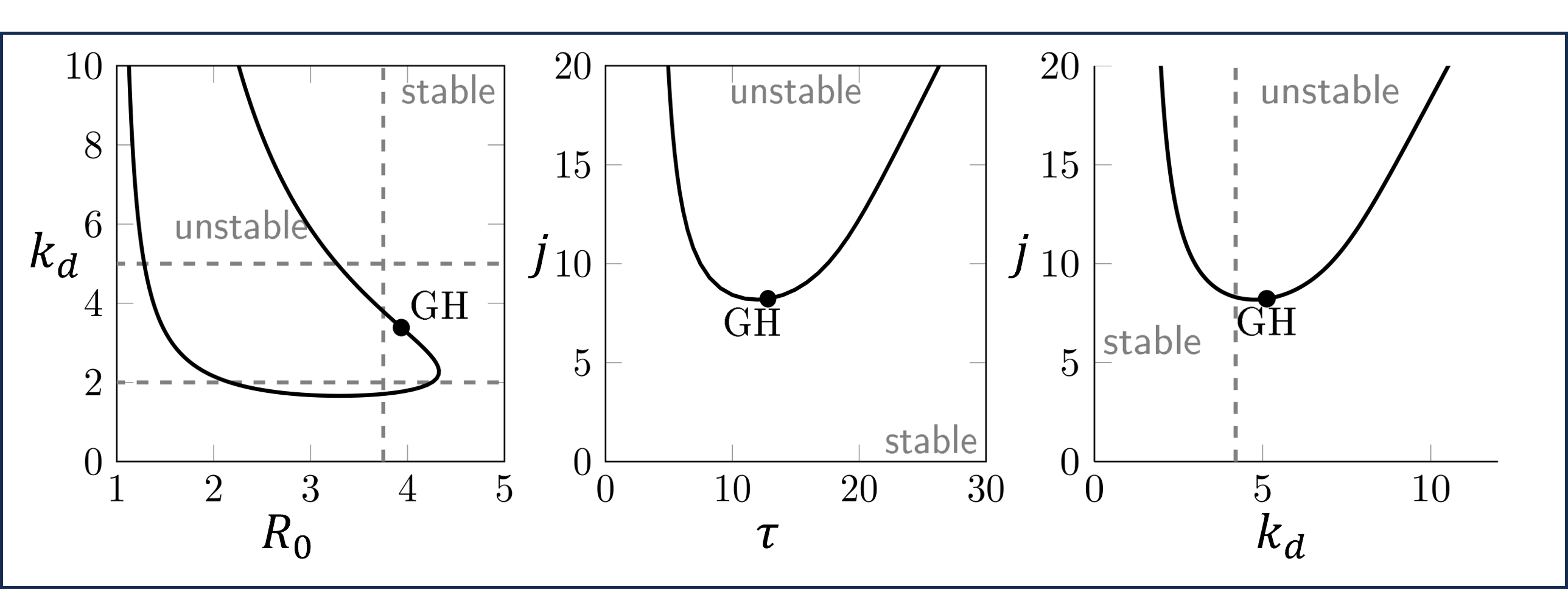}
    \caption{ \textbf{Stability of the EE in two-parameter planes} In all cases, the EE loses stability through a Hopf bifurcation. In Panel A,  $\tau=5$ and $j=20$ are fixed. In Panel B, we hold $R_0=2.2$ and $k_d=2$ fixed. In Panel C, we hold $R_0=2.2$ and $\tau=5$ fixed. The grey dashed lines correspond to parameter choices for the one-parameter continuations. }
    \label{fig:planes}
\end{figure}

We performed two-parameter bifurcation analyses of Eq.~\eqref{Eq:InfectedRenewalEquation} for generic gamma distributed durations of infection. Here, we consider $R_0 > 1,$ so that the DFE is unstable and the EE is strictly positive. We show the destabilization of the EE via a Hopf bifurcation and the resulting emergence of waves of infection in Figure~\ref{fig:planes}. In all cases, the boundary in two-parameter space separating the regions of stability and instability corresponds to a Hopf bifurcation. 
 
Panel A of Figure~\ref{fig:planes} illustrates the stability of the EE as the parameters representing the infectivity of the pathogen are varied. In particular, $k_d$ captures the waning infectivity of infectious individuals as they progress through infection. We note that, as the rate of this waning decreases and individuals are uniformly infectious through the infectious period, the system does not admit periodic solutions. Panel B illustrates how varying the duration of infection influences model dynamics. Here, for a fixed $\tau$, increasing the parameter $j$ acts to decrease the heterogeneity in the duration of infection with the gamma distribution approximates a delta distribution as $j \to \infty$. The resulting homogeneity in the duration of infection acts to synchronize the replenishment of the susceptible compartment, which drives waves of infections. Finally,  we vary the heterogeneity in both the duration of infection, via $j$, and the infectiousness, via $k_d$, while holding the mean duration of infection, $\tau$, and the infectivity, $R_0$, constant in Panel C. Here, we illustrate the role of heterogeneous infectivity and show that, even for the same mean duration of infection, different distributions of of the duration of infection lead to different epidemic dynamics. For example, if individuals are uniformly infectious with $k_d < 1$, then periodic dynamics can only occur due to the synchronization of the replenishment of the susceptible compartment. In Panels B and C, we vary the parameter $j$ continuously, which is not possible in the standard ODE formulation, and illustrates the additional utility of the pseudospectral approach.
 
 In all cases, a Generalized Hopf (GH) bifurcation separates the branch of supercritical Hopf bifurcations, characterised by the emergence of a stable branch of periodic solutions with corresponding waves of infection, from the subcritical Hopf bifurcations, characterised by the emergence of an unstable branch. These GH points are typically associated with the existence of a parameter region in which a stable equilibrium and a stable periodic solution coexist, so that the long-term dynamics of the system depend on the initial condition. 

\subsection{Genuine bistability} 

\begin{figure}[htp]
    \centering
    \includegraphics[trim= 4 10 5 10,clip,width=1\textwidth]{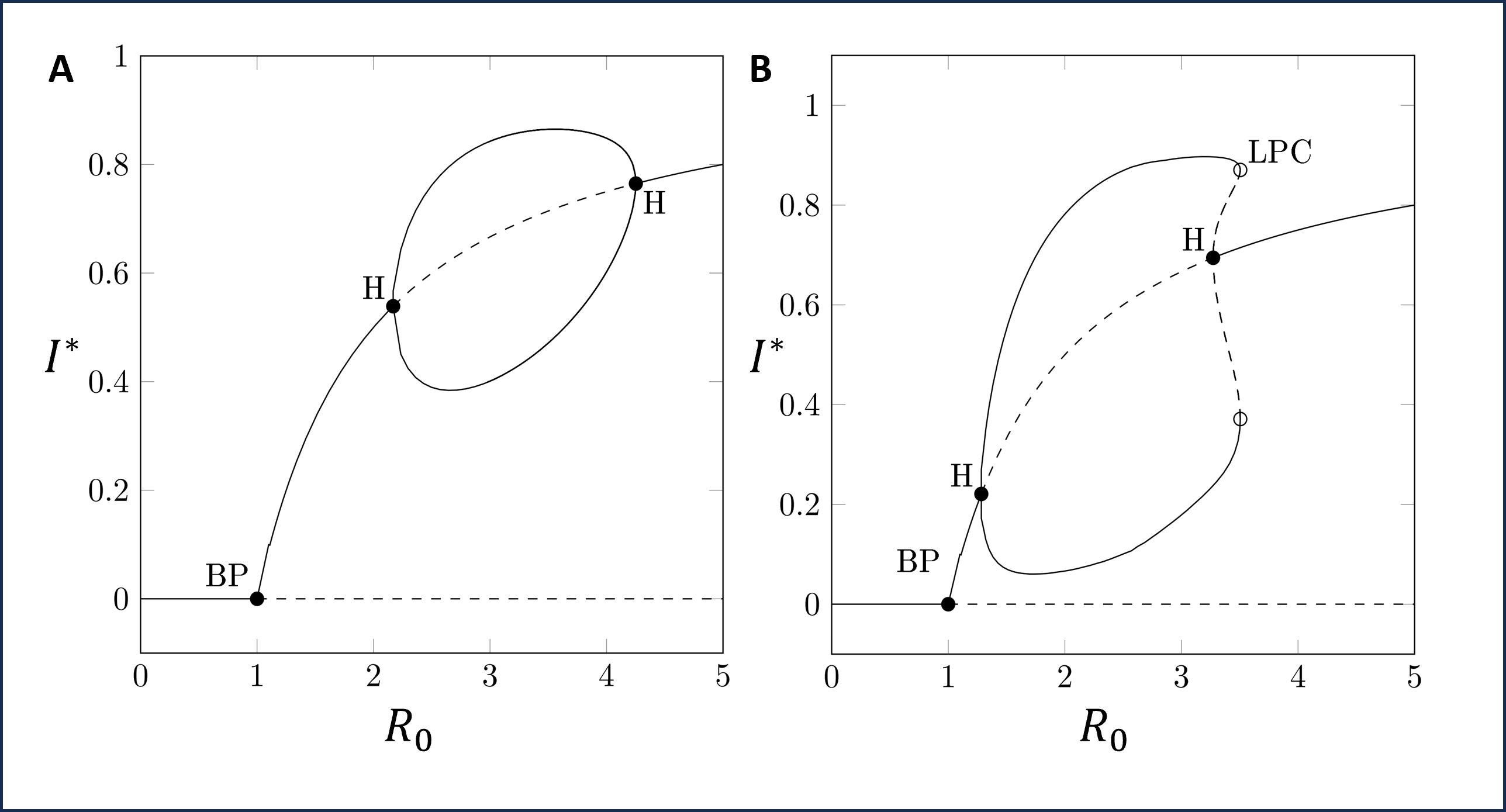}
    \caption{\textbf{One-parameter bifurcation diagram of the EE $I^*$ as a function of $R_0$} In both cases, the EE becomes stable following a transcritical bifurcation of the DFE at $R_0 =1$ and loses stability through a supercritical Hopf bifurcation. In Panel B, the EE undergoes a subcritical Hopf bifurcation which leads to the existence of an unstable periodic orbit and a region of bistability. In Panels A and B, we consider $\tau=5$ and $j=20$, while $k_d=2$ and Panel A and  $k_d=5$ in Panel B. These parameter choices correspond to the Panel A of Figure~\ref{fig:planes}. 
    }
    \label{fig:R_1p}
\end{figure}

Counterintuitively, rather than simply increasing the equilibrium value of $I$, increasing $R_0$ can either a stabilize or destabilize the EE. We illustrate this effect via a single-parameter bifurcation diagram with $R_0$ as the bifurcation parameter and plotting the size of the infected compartment via Eq.~\eqref{Eq:totalInfected} in Figure~\ref{fig:R_1p}. Here, we see the standard forward transcritical bifurcation of the DFE at $R_0 = 1$. However, as $R_0$ increases, the EE first loses stability, then becomes stable again via Hopf bifurcations. Here, a closed branch of periodic solutions exists between the two Hopf bifurcations. Panel~A shows a supercritical Hopf bifurcation where the resulting periodic solutions are always stable, while panel~B shows a subcritical bifurcation and the branch of periodic solutions is unstable near a Hopf point. This subcritical Hopf bifurcation results in a corresponding region of bistability, where the long-term dynamics of the epidemic depend on the initial conditions. This rich dynamical behaviour, of both periodic orbits and bistability between the EE and periodic orbits is not possible in the simple SIS model Eq.~\eqref{Eq:SISModelODE}. However, by including a gamma distributed duration of infection and waning infectivity, we observe the emergence of bistability in this simple model. 
 
\begin{figure}[htp]
    \centering
    \includegraphics[trim= 4 7 5 10,clip,width=1\textwidth]{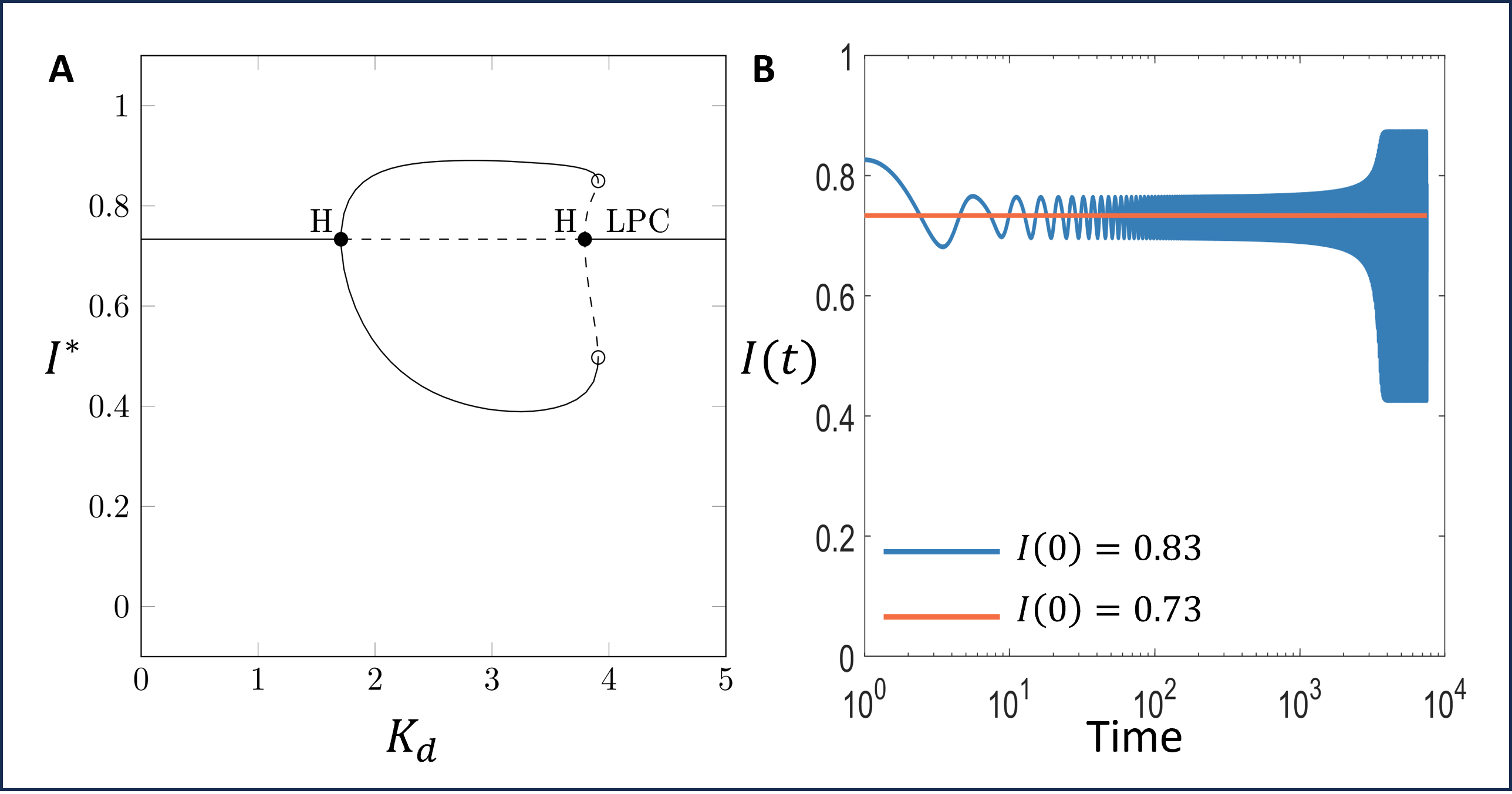}
    \caption{ \textbf{Genuine bistability in the simple epidemic model.}  In both Panel A and B, we take $R_0 =3.75, \tau=5$, and $j=20$. Panel A shows a one-parameter bifurcation diagram with $k_d$ as bifurcation parameter. The EE loses stability through a supercritical Hopf bifurcation which leads to a stable periodic orbit. This stable periodic orbit loses stability in a limit point of cycles bifurcation which results in a region of bistability of the EE and the stable periodic orbit. Panel B illustrates this bistability via time simulations of the simple epidemic model from two different initial conditions. Here, as $j = 20 \in \mathbb{N},$ the Erlang approximation is exact. }
    \label{fig:kd-bistability}
\end{figure}

We further demonstrate this bistability in Figure~\ref{fig:kd-bistability}, where the left panel shows a one-parameter bifurcation diagram with respect to $k_d$. We keep $R_0$ constant, so the number of infected individuals at the EE remains constant for all values of $k_d$, and take $j =20$. This EE undergoes stability switches (stable-unstable-stable) via two Hopf bifurcations. The first Hopf bifurcation is supercritical and results in the existence of stable periodic orbits while the second Hopf bifurcation is subcritical. This subcritical Hopf bifurcation delimits an interval of values of $k_d$ for which the stable EE coexists with a stable periodic solution. As mentioned, the total transmissibility $R_0$ and the rate of return to susceptibility are held constant, so the stability switches illustrated in Figure~\ref{fig:kd-bistability} are determined by the change in the infectiousness decay rate only. 

 We performed time simulations of the epidemic model with $k_d= 3.795$ to demonstrate the resulting dependence of long-term behaviour on initial conditions. Now, we recall that $j = 20,$ so the duration of infectiousness is Erlang distributed. Consequently, the ODEs in Eq.~\eqref{Eq:ErlangApproximationInfectedDistributeDDE} and Eq.~\eqref{Eq:LambdaErlangDE} are mathematically equivalent to Eq.~\eqref{Eq:InfectedDistributeDDE} and Eq.~\eqref{Eq:LambdaGenericDE}, respectively. In Figure~\ref{fig:kd-bistability}B, we show the time evolution from two distinct initial conditions
\begin{align*}
    I_1(0) = 1-\frac{1}{R_0+0.001} = 0.73, \quad \textrm{and} \quad I_2(0) = 1-\frac{1}{R_0+2} = 0.83.
\end{align*}
Here, we observe bistability from both the time simulation and bifurcation structure of the epidemic model. This bistability is a \textit{genuine} feature of the model that reflects how the inclusion of realistic waning infectiousness can drive complex dynamics in a simple epidemic model. 

\subsection{Spurious bistability}

\begin{figure}[htpb!]
    \centering
    \includegraphics[trim= 4 8 5 10,clip,width=1\textwidth]{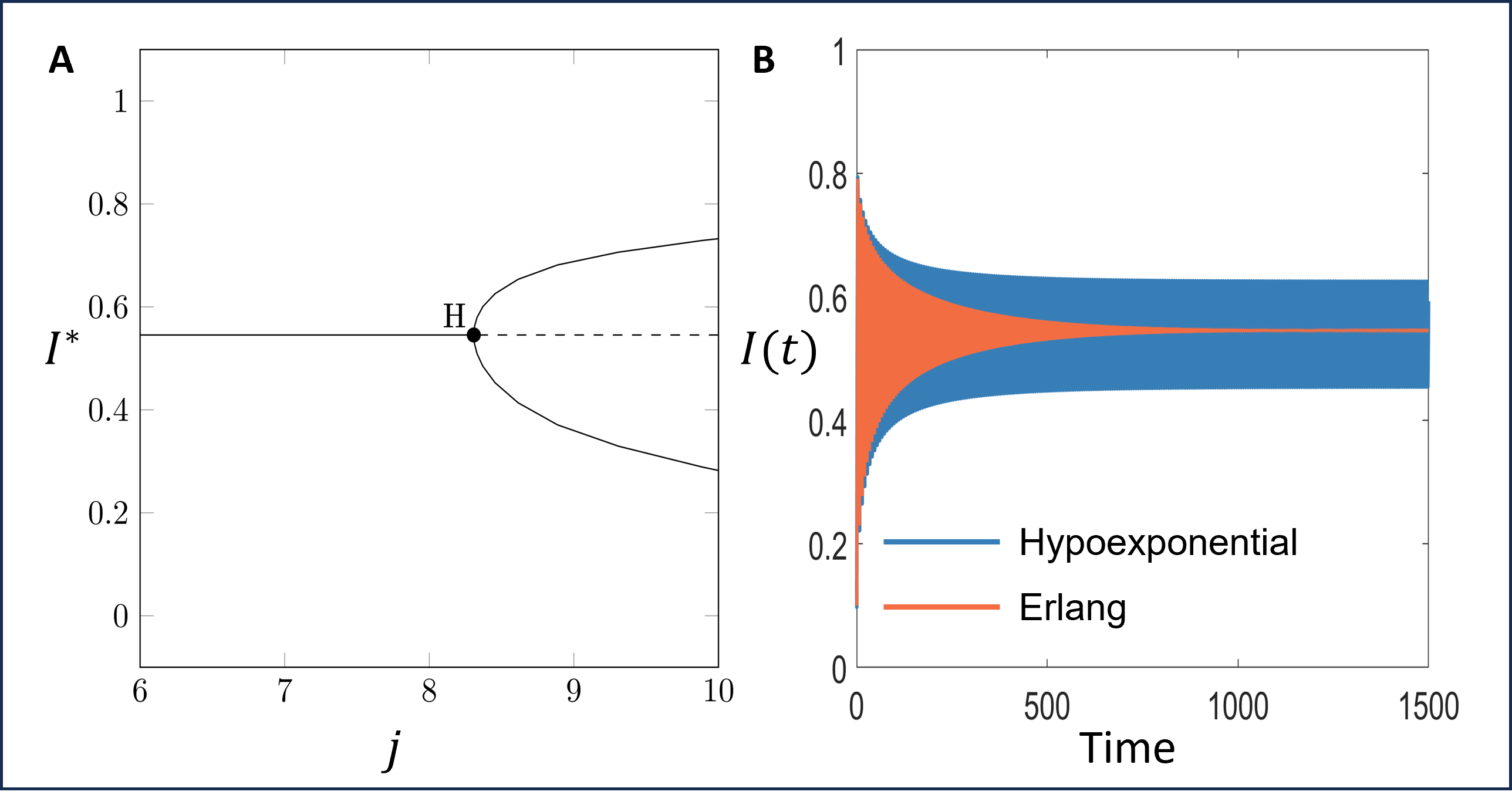}
    \caption{ \textbf{Spurious bistability in the simple epidemic model.}  In both Panel A and B, we take $R_0 = 2.2, \tau=5$, and $k_d = 4.2$. Panel A shows a one-parameter bifurcation diagram with $j$ as bifurcation parameter. The EE loses stability through a supercritical Hopf bifurcation which leads to a stable periodic orbit at $ j^* = 8.31$. Panel B shows two time simulations of this simple epidemic model using the Erlang and hypoexponential approximation to approximate the nonlinear RE with  $j = 8.32$. The Erlang approximation predicts spurious convergence to the EE while the hypoexponential approximation accurately captures the stable periodic orbit. }
    \label{fig:shape-one-parameter}
\end{figure}

We have thus far shown that the inclusion of waning infectivity can drive bistability in a simple epidemic model. However, in Figure~\ref{fig:planes}, we also demonstrated that the shape parameter, $j$, and thus the variance of the duration of infectivity, can impact the dynamical behaviour of the system. However, treating $j$ as a continuous parameter is not possible when approximating the underlying gamma distribution via the Erlang approximation and linear chain trick. Rather, treating $j$ as a continuous parameter is only possible when using the pseudospectral approach, which studies the underlying RE explicitly. We highlight the role of $j$ on the dynamics of this simple epidemic model in Figure~\ref{fig:shape-one-parameter}. Here, we identify a supercritical Hopf bifurcation at $j^* = 8.31$ and $\tau = 5$ using the pseudospectral approach. 

We next perform time simulations of the model using both the Erlang and Hypoexponential ODE approximations in Section~\ref{Sec:ODEApproximations} to illustrate the potential effects of these approximations for $j = 8.32$. Here, $j > j^*$, so our bifurcation analysis indicates the existence of a stable periodic orbit. We perform the hypoexponential approximation in Eq.~\eqref{Eq:InfectedHypoODE} and Eq.~\eqref{Eq:LambdaHypoODE} using the equivalent value of $\sigma^2 = \tau^2/j.$ Conversely, to perform the Erlang approximation, we round the shape parameter $j$ to the nearest integer. However, this rounding moves the Erlang approximation across the Hopf bifurcation and the time simulations predict convergence to the EE. Conversely, the hypoexponential approximation, obtained by taking $\tau = 5$ and $\sigma^2 = \tau^2/j$, captures the correct periodic behaviour. Consequently, the time dynamics of the simple epidemic model depend on which approximation is used, leading to \textit{spurious} bistability. Here, it is important to note that the time simulations are otherwise identical and this divergence in long term behaviour is solely due to the approximation error. 

\section{Discussion}

We developed a simple extension of the standard SIS ODE model that includes time-since-infection dependence in both the infectivity and the return to susceptibility rates. By decoupling the infectivity rate from the rate at which infected individuals return to the susceptible class, our SIS model can also describe SIRS-type situations in which infected individuals experience a period in which they are not infectious anymore, but not yet susceptible. 

Through numerical bifurcation analysis, we have shown that the inclusion of these two mechanisms can lead to the destabilization of the EE via a Hopf bifurcation and the emergence of periodic solutions. Further, we demonstrate the existence of bistability between stable periodic orbits and the EE following a subcritical Hopf bifurcation. 
To our knowledge, our results are the first example of bistability arising in a simple SIS model with no other mechanisms like seasonal forcing or boosting of immunity. Consequently, our findings further highlight how temporary immunity at the individual level can alone drive complex dynamics at the population scale. 

Our results highlight how the distribution of the duration of infection impacts model dynamics. Typically, the mean duration of infection is used to derive the rate at which individuals return to the susceptible pool. However, here, we showed that the higher order moments of this distribution, such as the variance, can lead to instability of the EE and the emergence of a stable periodic orbit, and to bistability. For example, we showed that if the variance is small, the resulting homogeneity in the duration of infection can lead to oscillatory dynamics. These results are particularly relevant in the context of compartmental epidemic models, which correspond to implicitly assuming an Erlang distributed infectious period. Indeed, we showed that implementing the common Erlang approximation can alter the effective location of Hopf bifurcations. We showed that this approximation can lead to the existence of spurious bistability, wherein the ODE approximation incorrectly predicts the existence of a stable EE. We next demonstrated how the hypoexponential approximation \citep{Cassidy2022}, which matches both the first and second moments of the underlying distribution, does not exhibit this spurious bistability while retaining the numerical tractability of ODE models. 

We are able to continuously vary the variance of the duration of infection by using a pseudospectral approximation. Specifically, this numerical bifurcation approach facilitates the direct investigation of the underlying RE describing the epidemic without making additional, simplifying assumptions. Taken together, our approach illustrates how modellers can derive simple epidemic models from either the cohort or generational perspective before directly analysing these models using computationally efficient techniques. These pseudospectral techniques are increasingly used in the analysis of dynamics and bifurcations of physiologically structured models, and are applicable to both the RE and structured PDE setting \citep{Scarabel2021Vietnam, DiekmannScarabelSize, DiekmannScarabelAge,Breda2016EJQTDE}. 

Here, we have considered a deliberately simple epidemic model to illustrate how the inclusion of variable infectivity and waning immunity can drive complex dynamics in a SIS model. However, while this model is sufficient to illustrate these dynamics, it is likely too simple to be applicable in many data-driven settings. For example, the function $\beta(a)$ will likely have compact support in most data-driven applications, rather than the simple example we considered here. Further, natural extensions of our model include incorporating waning vaccine-induced immunity and relaxing the assumption of a well-mixed population. Despite these limitations, our simple model is sufficient to illustrate how decoupling infectivity from the duration of infection can lead to the emergence of bistability in a simple SIS model. 

More importantly, our results illustrate how approximating the structured model using numerically tractable, compartmental ODE models can lead to incorrect predictions regarding the asymptotic behaviour of the model, and suggest that caution is needed when using these compartmental epidemic models. Consequently, we provide two complementary approaches to avoid these spurious predictions, by either studying the bifurcation structure of the model directly using the pseudospectral approach or by using the hypoexponential approximation to perform time simulations. Consequently, we illustrate how to include the biologically relevant mechanisms independently driving infectivity and waning immunity within a simple SIS modelling framework while retaining the ability to simulate and analyse the resulting model. 

\section*{Acknowledgements}
TC and FS were partially supported by the Engineering and Physical Sciences Research Council via the Mathematical Sciences Small Grant UKRI170: ``The dynamics of waning and boosting of immunity: new modelling and numerical tools''. FS is a member of the Computational Dynamics Laboratory (CDLab, University of Udine), of INdAM research group GNCS, of UMI research group Mo\-del\-li\-sti\-ca Socio-Epidemiologica, and of JUNIPER (Joint UNIversities Pandemic and Epidemiological Research). HC was supported by a University of Leeds Summer Bursary.

\clearpage

\end{document}